\documentclass[smallextended]{svjour3}       

\smartqed  
\usepackage{mathptmx}      
\usepackage{latexsym}
\usepackage{graphicx}
\usepackage[centertags]{amsmath}
\usepackage{amssymb}
\usepackage{amsfonts}
\usepackage{soul}
\usepackage{multicol}
\usepackage{mathrsfs}

\newtheorem{propo}[theorem]{Proposition}

\newcommand{\bt}{\begin{theorem}}
\newcommand{\et}{\end{theorem}}
\newcommand{\bl}{\begin{lemma}}
\newcommand{\el}{\end{lemma}}
\newcommand{\bp}{\begin{proposition}}
\newcommand{\ep}{\end{proposition}}
\newcommand{\bc}{\begin{Corollary}}
\newcommand{\ec}{\end{Corollary}}
\newcommand{\bdefn}{\begin{Definition}}
\newcommand{\edefn}{\end{Definition}}
\newcommand{\br}{\begin{Remark}}
\newcommand{\er}{\end{Remark}}
\newcommand{\be}{\begin{equation}}
\newcommand{\ee}{\end{equation}}
\newcommand{\ba}{\begin{array}}
\newcommand{\ea}{\end{array}}
\newcommand{\bea}{\begin{eqnarray}}
\newcommand{\eea}{\end{eqnarray}}
\newcommand{\bean}{\begin{eqnarray*}}
\newcommand{\eean}{\end{eqnarray*}}
\newcommand{\beit}{\begin{itemize}}
\newcommand{\eeit}{\end{itemize}}
\newcommand{\ben}{\begin{enumerate}}
\newcommand{\een}{\end{enumerate}}

\catcode`\@=11
\def \theequation{\@arabic{\c@section}.\@arabic{\c@equation}}    
\def \thetheorem{\@arabic{\c@section}.\@arabic{\c@theorem}}
\def \theenumi{\@roman{\c@enumi}}
\catcode`\@=12

\newcommand{\ds}{\displaystyle}

\newcommand \noi{\noindent}

\newcommand{\C}{\mathcal{C}}

\newcommand{\intot}{\int_{\Omega_t}}

\newcommand{\nro}{\textordmasculine\ }

\newcommand{\lp}{\left(}
\newcommand{\rp}{\right)}

\newcommand{\lc}{\left[}
\newcommand{\rc}{\right]}
\newcommand{\lan}{\langle}
\newcommand{\ran}{\rangle}

\newcommand{\Om}{\Omega}

\newcommand{\intom}{\int_\Om}

\newcommand{\ov}{\overline}
\newcommand{\La}{\Delta}

\newcommand{\vphi}{\varphi}

\newcommand{\n}{\nabla}

\newcommand{\F}{\mathcal{F}}

\newcommand{\R}{\mathbb{R}}

\newcommand{\ra}{\rightarrow}

\begin{document}

\title{A shape optimization problem for the $p$-Laplacian
\thanks{We thank Comisi\'on Nacional de Investigaci\'on Cient\'ifica y
  Tecnolog\'ica(CONICYT) for financial support through the project
  FONDECYT Regular N\nro 1090305 as also National Board of Higher
  Mathematics (NBHM), India and Indian Institute of Science Education and
  Research (IISER), Pune.}}

\titlerunning{$p$-Laplacian optimization}        

\author{Anisa Chorwadwala         \and
        Rajesh Mahadevan 
}


\institute{Anisa Chorwadwala, \at
              Indian Institute of Science Education and Research (IISER),  \\
             Central Tower, Sai Trinity Building, \\
Sutarwadi Road, Pashan, Pune 411021, INDIA.\\
              \email{anisa@iiserpune.ac.in}           
           \and
           Rajesh Mahadevan, \at
              Departamento de Matem\'atica, \\
              Facultad de Cs. F\'isicas y Matem\'aticas,\\
              Av. Esteban Iturra s/n, Barrio Universitario,\\
              Casilla 160C, Concepci\'on, 
              CHIL\'E.\\
              \email{rmahadevan@udec.cl}           
}
\date{Received: date / Accepted: date}

\maketitle
\begin{abstract}
It is known that the torsional rigidity for a punctured ball, with the
puncture having the shape of a ball, is minimum when the balls are
concentric and the first eigenvalue for the Dirichlet Laplacian for
such domains is also a maximum in this case. These results have been
obtained by Ashbaugh and  Chatelain (private communication), Harrell
et. al. \cite{HKK}, Kesavan \cite{Kesavan} and Ramm and Shivakumar
\cite{Ramm-Shivakumar}. In this paper we extend these results to 
the case of $p$-Laplacian for $1<p<\infty$. For proving these results,
we follow the same line of ideas as in the aforementioned articles,
namely, study the sign of the shape derivative using the moving  plane
method and comparison principles. In the process, we obtain some
interesting new side results such as the Hadamard perturbation formula
for the torsional rigidity functional for the Dirichlet $p$-Laplacian,
the existence and uniqueness result for a nonlinear pde and some
extensions of known comparison results for nonlinear pdes.
\end{abstract}

\section{Introduction}\label{S:intro}
\setcounter{equation}{0}
\setcounter{theorem}{0}
The $p$-Laplacian $\Delta_p$ is the non-linear operator defined as 
$\Delta_p f = \mathrm{div}(|\nabla f|^{p-2}\nabla f)$. Let $B_1$ be an
open ball in $\R^N$. Let $B_0$ be another open ball whose closure is
contained in $B_1$, and is free to move inside $B_1$. Let $\Omega=B_1 \setminus \ov{B_0}$. 
We consider the following domain optimization problems:
\begin{enumerate}
\item Given $y \in W^{1,p}_0(\Omega)$, the unique solution of the equation  
\be\label{stationary}
\left.
\begin{aligned}
-\Delta_p u = & \; 1 \quad\text{in $\Omega$,} \\
        u = & \; 0 \quad\text{on $\partial\Omega$},
\end{aligned}~~~ \right\} 
\ee
we are interested in minimizing the $p$-torsional rigidity 
\be\label{trsnlrgdty}
E(\Om):=\intom |\nabla y|^p \, dx=\intom y \,dx
\ee
with respect to the position of the hole $B_0$. 

\item Given the eigenvalue problem
\be\label{evp}
\left.
\begin{aligned}
-\Delta_p u = & \; \lambda |u|^{p-2}\,u \quad\text{in $\Omega$,} \\
        u = & \; 0 \quad\text{on $\partial\Omega$}
\end{aligned}~~~ \right\} 
\ee
whose principal eigenvalue 
is
\be\label{1stevplap}
\lambda_1(\Omega) :=
\inf \left\{ \left.\dfrac{\|\nabla \vphi\|^p_{L^p(\Om)}}{\|\phi
  \|^p_{L^p(\Om)}}\;\right|\; \vphi \in W^{1,p}_0(\Om)\right\}\,,
\ee
we are interested in maximizing $\lambda_1(\Omega)$ with respect to
the position of $B_0$. 
\end{enumerate}
The following results were obtained, in the linear case, i.e., for
$p=2$, by Ashbaugh and Chatelain (private communication), Harrell
et. al. \cite{HKK}, Kesavan \cite{Kesavan}, Ramm and Shivakumar
\cite{Ramm-Shivakumar}: the torsional rigidity is minimum if and only
if $B_0$ and  $B_1$ are concentric. Also, the first eigenvalue
$\lambda_1$ of problem (\ref{1stevplap}) attains its maximum if and
only if the balls are concentric.   

The analogues of these results for manifolds were obtained in Anisa 
and Aithal \cite{Anisa-Aithal} in the setting of  space-forms (complete
simply connected Riemannian manifolds of constant sectional 
curvature) and in Anisa and Vemuri (On two functionals connected to
the Laplacian in a class of doubly connected domains in rank-one
symmetric spaces of non-compact type, preprint) in the setting
of rank-one symmetric spaces of non-compact type. We extend these
results, in a different direction, to the non-linear setting.  Our
main results are Theorem \ref{trsnlrgdtyrslt} and Theorem \ref{evprslt}.

The proofs in \cite{Kesavan,Ramm-Shivakumar} rely on {\em shape
  differentiation} \cite{Sokolowski-Zolesio}, the {\em moving plane
  method} \cite{Berestycki-Nirenberg,Gidas-Ni-Nirenberg} and various
maximum principles. In the non-linear case, carrying out this program
involves several technical difficulties. We develop the shape calculus
for the torsional rigidity function for $p$-Laplacian. A formula for
the Hadamard perturbation of the first Dirichlet eigenvalue for the
$p$-Laplacian is given. This, however, is not new and may also be seen
in the works of Garc\'ia Meli\'an and Sabina de Lis \cite{GMSlS}, 
Lamberti \cite{L} and Ly \cite{Ly}. For the Steklov eigenvalue this is
done in Del Pezzo and Fern\'andez Bonder \cite{dPFB}. Subsequently, we
analyze the sign of the shape derivative. We do this by proving a
suitable strong comparison result. In the case of the eigenvalue
problem, before this, we also need to prove a general weak comparison
principle for the $p$-Laplacian with non-vanishing boundary condition
(cf. Theorem \ref{WCPbetainc}).  This result is new and can be of
independent interest in itself. An existence and uniqueness result for
a nonlinear pde is required for applying this comparison principle and
this result is also proved (cf. Proposition \ref{cdnA-2}).  

The Section  \ref{prel} establishes notations, contains some
definitions and technical preliminaries. In Section
\ref{comparison-plap}, we recall some existing weak and strong
comparison principles for the $p$-Laplacian and prove an extension of
a weak comparison principle. In Section \ref{anlpde}, we prove the
existence and uniqueness of non-negative solution for a nonlinear pde  
needed for an application of the comparison principle.  In  Section
\ref{shpdrvtvs}, following \cite{Sokolowski-Zolesio}  we obtain the
Hadamard perturbation formula for the torsional rigidity functional
\eqref{trsnlrgdty} and for the first eigenvalue  of the Dirichlet
$p$-Laplacian (\ref{1stevplap}). Finally, in  Section \ref{mainrslt}
we prove the  main results by analyzing the sign of the shape
derivatives.         

\section{Prelminaries}\label{prel}
\setcounter{equation}{0}
\setcounter{theorem}{0}
In this section we introduce some definitions and recall some results
which will be used later on. 

\bigskip
\noindent
{\sc Shape Derivative: } Given a functional $J$ which depends on the
domain $\Omega$ (usually, a smooth open set in $\mathbb{R}^N$) and
given, a variation of the domain $\Omega$ by a fairly smooth
perturbative vector field $V$ which has its support in a
neighborhood of $\partial \Omega$, the infinitesimal 
variation of $J$ in the direction $V$ is defined as 
\be
\label{shapeder} 
\ds{J'(\Omega;V)=\lim_{t\to 0} \frac{J( \Omega_t ) -
    J(\Omega)}{t}\,} 
\ee
where $\Omega_t$ is the diffeomorphic image $\Phi_t(\Omega)$ of
$\Omega$ under the smooth perturbation of identity $\Phi_t(x)= (I + t
V)(x)$. 

The shape derivative is a tool widely used in problems of optimization
with respect to the domain as it permits to understand the
variations of shape functionals (cf. Simon \cite{S},
\cite{Sokolowski-Zolesio}).  

We define : $B(t) := (D\Phi_t)^{-1}$, $\gamma(t):= |\mathrm{det}
D\Phi_t|$ and $A(t):= \gamma(t)\, B(t)B(t)^*$ where $B(t)^*$ shall 
denote the transpose of $B(t)$. It will be convenient to denote
$\gamma(t)$, $B(t)$, $B(t)^*$ and $A(t)$ respectively, by $\gamma_t$,
$B_t$, $B_t^*$ and $A_t$. We observe that  
\be
\label{gradients} D\Phi_t = I + t \, DV \, \qquad \qquad \ (D\Phi_t)^* 
= I +  t \ (DV)^*\\  
\ee
and so, $B_t$, $B_t^*$, $A_t$, $\gamma_t$ and $F_t$ are analytic
functions of $t$ near $t=0$. We record that 
\bea
\label{drvtvgammaat0} \gamma'(0) & = & \mathrm{div} \, V\\
\label{drvtvBt0} (B_t^*)'(0) & = &- (DV)^*.
\eea
So, for small $t$, we have
\bea
\label{asygammat} \gamma_t & \approx & \gamma(0) + t \gamma'(0) = 1 +
t \, \mathrm{div} V 
\eea
Also, for $t$ sufficiently small say $|t| < t_0$, there exists a
constant $C  >0$ such that 
\be\label{bdphit}
|(D\Phi_t)^* \xi |  \le C  |\xi| \mbox{ for all } \xi
\in \R^n\,. 
\ee
Consequently, by substituting $B_t^* \eta$ for $\xi$, $\eta$
arbitrary in $\R^n$, we have 
\be\label{bdBt}
|B_t^* \eta |  \ge C^{-1}  |\eta| \mbox{ for all } \eta \in \R^n\,.
\ee

\smallskip
\noindent
{\sc Pucci-Serrin identity: } We shall find it very useful to employ
the extended version of the Pucci-Serrin identity proved by
Degiovanni et. al. \cite{DMS} which gives the following identity for
the $p$-Laplacian. Assume that $u \in C^1(\overline{\Omega})$ is a
solution of the equation
\be\label{plapeq}
\left.
\begin{aligned}
-\Delta_p u  & =  f ~& \mbox{ in } & ~\Om,  \\
u & =  0 ~& \mbox{ on } & ~\partial \Om. 
\end{aligned} 
~~~ \right\} 
\ee
Then, for all $V \in C^1(\overline{\Omega})$ the following identity
holds  
\be \label{PucciSerrinC1} 
\ds{-\dfrac{(p-1)}{p} \intom \left|\dfrac{\partial
    u}{\partial \nu}\right|^p  \, V \cdot n \ dS  = \intom \mathrm{div}
  V  \,  \dfrac{|\n u|^p}{p} \ dx -\intom \lan (DV)^* \n u \,,\, |\n
  u|^{p-2} \n u \ran \ dx+ \intom V \cdot \n u \, f \ dx  }
\ee
The Pucci-Serrin identity may be obtained  by using $V \cdot \n u$ as
a test function in \eqref{plapeq} and after several integration by
parts whenever $u \in  C^1(\overline{\Omega}) \cap
C^2(\Omega)$. However, by standard regularity results for solutions of
the $p$-Laplacian equation, they are known to belong to only
$C^{1,\alpha}(\overline{\Omega})$ (cf. Tolksdorff \cite{Tolk2}) as the
coefficients $|\n u|^{p-2}$ degenerates near the critical points of
$u$.  This formula can be justified by regularizing the coefficient
first and then passing to the limit cf. \cite{DMS} (see also
Garc\'ia-Melian  and Sabina de Lis \cite{GMSlS} and the work of Del
Pezzo and Fern\'andez Bonder \cite{dPFB} for such arguments). 

\smallskip
\noindent
{\sc A positive definite matrix:} Define a strictly
convex function $\Gamma : \R^N \ra \R$ by $\Gamma(x)=
\frac{|x|^p}{p}$. Let $A= D\Gamma$. Then $A= (A_1, A_2, \ldots, A_N) :
\R^N \ra \R^N$ and is given by 
\be\label{2.10}
A(x)= |x|^{p-2}x.
\ee 
Clearly, $A \in \C^\infty(\R^N \setminus \{0\})$. The matrix
$\mathcal{A}:=\left[ \dfrac{\partial A_i}{\partial
    x_j}(x)\right]_{i,j =1}^N$ corresponds to the symmetric matrix
$|x|^{p-2} \mbox{Id} + (p-2) |x|^{p-4} x \otimes x$, which is the
Hessian of the convex function $\Gamma$. It can be seen that $(p-1)|x|^{p-2}$ and $|x|^{p-2}$ are eigenvalues
of $\mathcal{A}$ with multiplicity one and $(n-1)$
respectively. Therefore, for any $\xi \in \R^n$, we have  
\be \label{2.11}
<\mathcal{A}\xi,\xi> \geq \min\{1,
 p-1\} |x|^{p-2} |\xi|^2.
\ee
\section{Comparison Theorems for the
  $p$-Laplacian}\label{comparison-plap} 
\setcounter{theorem}{0}
\setcounter{equation}{0}
Let $\Om \subset \R^N$ be a bounded domain with smooth boundary. Let
$\beta : \Omega \times \R \rightarrow \R$ be continuous function and
$u \longmapsto \beta(x,u)$ is locally Lipschitz on $\R \setminus \{0\}$
uniformly for $x \in \Om$ and assume that $\dfrac{\partial
\beta}{\partial u}$ is of constant sign for all $(x,u) \in \Om \times
\lp\R \setminus \{0\}\rp$. Let $f,g \in W^{-1, \frac{p}{p-1}}(\Om)$, 
$f^\prime,g^\prime \in W^{1-\frac{1}{p},p}(\partial \Om)$ with $f \geq
g$ in $\Om$ (in the sense of distributions), $f^\prime \geq g^\prime$
on $\partial \Om$. Let $u,v \in W^{1,p}(\Om)$ solve (in the weak sense) 
\begin{equation}\label{nlpdes}
\begin{aligned}
-\La_p u &= \beta(x,u) + f(x), & -\La_p v &= \beta(x,v) + g(x) &
\mbox{ in } & \Om,\\ 
u&= f^\prime, & v &=g^\prime  & \mbox{ on } & \partial \Om. 
\end{aligned}
\end{equation}
Then one is interested in the following comparison results: 
\begin{itemize}
\item[(WCP)] Weak Comparison Principle: Is it true that $u \geq v$ in
  $\Om$? 
\item[(SCP)\,] Strong Comparison Principle: If $u,v \in
  \C^1(\ov{\Om})$, $u \not \equiv v$, $u \geq v$ in $\Om$, is it true
  that $u > v$ in $\Om$ and $\frac{\partial u}{\partial n}(x_0)<
  \frac{\partial v}{\partial n}(x_0)$ for any $x_0 \in \partial \Om$?
  Here, $n$ is the unit outward normal to $\Om$ on $\partial \Om$. 
\end{itemize}
The Weak Comparison Principle (WCP) holds when $\dfrac{\partial
  \beta}{\partial u} \le 0$ for which we refer to Tolksdorff
\cite{Tolksdorff}.  

\smallskip
\noi
The Weak Comparison Principle also holds when $\dfrac{\partial
  \beta}{\partial u} \geq 0$ under the following assumptions and for
Dirichlet boundary data: 
\begin{itemize}
\item[(A-1)] $\dfrac{\partial \beta}{\partial u} \geq 0\; \forall
  \;(x,u) \in \Om \times \lp\R \setminus \{0\}\rp$, $\beta(x,0) \geq 0
  \; \forall \; x \in \Om$. 
\item[(A-2)] The problem 
\begin{equation*}
\begin{aligned}
-\La_p u &= \beta (x,u) + f & \mbox{ in } & \Om,\\
u & = 0 & \mbox{ on } & \partial \Om.
\end{aligned}
\end{equation*}
(where $f \in L^\infty(\Om)$, $f \geq 0$ in $\Om$) admits a unique
non-negative solution $u \in W_0^{1,p}(\Om)$. 
\item[(A-3)] $f,g \in L^\infty(\Om)$, $0 \leq g \leq f$ on $\Om$ and
  $0 = g^\prime = f^\prime$ on $\partial \Om$.
\end{itemize}
This result is proved in \cite{Cuesta}. However, for our purposes the
zero Dirichlet data assumption in \newline (A-2) and (A-3) is too restrictive. We
show that this result also holds for inhomogeneous Dirichlet boundary
data, that is, by relaxing the condition (A-2) and (A-3) to (A-2') and
(A-3') respectively:
\begin{itemize}
\item[(A-2')] The problem 
\begin{equation*}
\begin{aligned}
-\La_p u &= \beta (x,u) + f & \mbox{ in } & \Om,\\
u & = f' & \mbox{ on } & \partial \Om.
\end{aligned}
\end{equation*}
(where $f \in L^\infty(\Om)$, $f \geq 0$ in $\Om$ and $f'\ge 0$ on
$\partial \Omega$) admits a unique non-negative solution $u \in
W^{1,p}(\Om)$.  
\item[(A-3')] $f,g \in L^\infty(\Om)$, $0 \leq g \leq f$ on $\Om$ and
  $0 \le g^\prime \le f^\prime$ on $\partial \Om$.\\
\end{itemize}
We prove the following results along the same lines as in
\cite{Cuesta}. 
\begin{theorem}\label{WCPbetainc}
Let the assumptions (A-1) (A-2') and (A-3') hold then the WCP holds for
bounded solutions.
\end{theorem}
\noindent {\it Proof.} Let us denote $L_+^\infty(\Om)= \left\{ h \in
L^\infty(\Om)\,|\, h \geq 0 \mbox{ in } \Om\right\}$. Given $f \in
L_+^\infty(\Om)$ and $f^\prime \in W^{1-\frac{1}{p},p}(\partial \Om)$
with $f' \ge 0$ on $\partial \Om$, define the nonlinear operator
$T_{f,f^\prime}$  on $L_+^\infty(\Om)$ by letting
$T_{f,f^\prime}(u)=v$, where $v$ is the 
weak solution of  
\begin{equation}\label{defTffprime}
\begin{aligned}
-\La_p v &= \beta (x,u) + f & \mbox{ in } & \Om,\\
v & = f^\prime & \mbox{ on } & \partial \Om
\end{aligned}
\end{equation}
Since $\dfrac{\partial \beta}{\partial u} \geq 0$ and $u \ge 0$, it
follows that $\beta(x,u) \geq \beta(x,0)\geq 0$. So, the right hand
side in \eqref{defTffprime} is non-negative as also the boundary
data. By appealing to the WCP proved by Tolksdorff \cite{Tolksdorff}
we conclude that indeed $T_{f,f^\prime}(u)=v \ge 0$ and
$T_{f,f^\prime}$ maps $L_+^\infty(\Om)$ into itself.\\ 
\noindent {\it Claim.} Let $f_1, f_2, u_1,u_2 \in L_+^\infty(\Om)$. If 
$f_1 \leq f_2$, $u_1 \leq u_2$ and $f_1^\prime \leq f_2^\prime$ then
$T_{f_1,f_1^\prime}(u_1) \leq T_{f_2,f_2^\prime}(u_2)$\\ 
Indeed, following the condition  $\dfrac{\partial \beta}{\partial u}
\geq 0$ we conclude that  $f_1^* := \beta (x,u_1) + 
f_1 \leq \beta(x,u_2) + f_2 =: f_2^*$. Let
$v_1=T_{f_1,f_1^\prime}(u_1)$ and $v_2=T_{f_2,f_2^\prime}(u_2)$. Then 
\begin{equation*}
\begin{aligned}
-\La_p v_1 &=  f_1^*, \quad \quad &  -\La_p v_2 &=  f_2^* &~~~\quad
\mbox{ in } & \Om,\\ 
v_1 & = f_1^\prime, \quad \quad & v_2 & = f_2^\prime &~~~ \quad
\mbox{ on } & \partial \Om. 
\end{aligned}
\end{equation*}
So, again by the weak comparison result proved in \cite{Tolksdorff} 
we obtain $v_1 \leq v_2$ in $\Om$. This proves the claim. \\ 
Now, let $u,v$ be bounded solutions of the non-linear pdes in
\eqref{nlpdes}. To begin with, $T_{f,f^\prime} (u)=u$ and
$T_{g,g^\prime} (v)=v$. Now, using the claim we obtain the
inequalities,  
\[
0 \leq T_{f,f^\prime} (0) \leq T_{f,f^\prime} (u)=u\,, \qquad \qquad 0
\leq T_{g,g^\prime} (0) \leq T_{g,g^\prime} (v)=v\, 
\]
We can then show by an inductive application of the claim that
following chains of inequalities hold  
\bea
\label{chainu} &&0 \leq T_{f,f^\prime} (0) \leq T_{f,f^\prime}^2 (0)
\leq \cdots \leq T_{f,f^\prime}^n (0) \leq \cdots \leq u
=T_{f,f^\prime} (u) \\
\label{chainv} &&0 \leq T_{g,g^\prime} (0) \leq T_{g,g^\prime}^2 (0)
\leq \cdots \leq T_{g,g^\prime}^n (0) \leq \cdots \leq v
=T_{g,g^\prime} (v)
\eea
The pointwise limits $\ds u^*(x)= \lim_{n\longrightarrow \infty}\lc
T_{f,f^\prime}^n (0)\rc(x)$ and $\ds v^*(x)= \lim_{n\longrightarrow
  \infty}\lc T_{g,g^\prime}^n (0)\rc(x)$ exist and must clearly
satisfy $T_{f,f^\prime} (u^*)=u^*$ and $T_{g,g^\prime} (v^*)=v^*$
respectively. So, by the uniqueness assumption in (A-2'), it follows
that  $u^*=u$ and $v^*=v$.

Again, by applying the claim above, for any $n \geq 0$, we obtain 
$T_{g,g^\prime}^n (0)\leq T_{f,f^\prime}^n (0)$. Therefore, upon
taking the limit as $n$ goes to infinity we obtain $v \leq u$. This
proves the theorem. \hfill $\square$
\section{Existence and uniqueness for a nonlinear Dirichlet problem
}\label{anlpde} 
\setcounter{theorem}{0}
\setcounter{equation}{0}
Let $\lambda_1$ be the first eigenvalue of the Dirichlet $p$-Laplacian
as in \eqref{1stevplap} on a bounded domain $\Omega$. Let $\mathcal{O}$ be an
open proper subset of $\Omega$. We prove the existence and uniqueness
result for a nonlinear partial differential equation on $\mathcal{O}$
given Dirichlet data $f^\prime \ge 0$ on $\partial \Omega$. This shall
be needed for applying the comparison principle of the previous
section, later in Section \ref{mainrslt}.   
\begin{propo}\label{cdnA-2}
Given $f^\prime \in W^{1-\frac{1}{p},p}(\partial \mathcal{O})$ and $f^\prime
\ge 0$ on $\partial \mathcal{O}$, the problem 
\be \label{nlpde}
\left.
\begin{aligned}
-\La_p w &= \lambda_1\, |w|^{p-2}\,w  & \mbox{ in } & \mathcal{O},\\
w & = f^\prime & \mbox{ on } & \partial \mathcal{O}.
\end{aligned}~~~ \right\} 
\ee
admits a unique non-negative solution. 
\end{propo}
\noindent {\it Proof.} \quad Let us first prove that if a solution
exists then it is non-negative. Let $u$ be a solution of the above
problem. As $u \ge 0$ on $\partial \mathcal{O}$, we obtain that $u^- \in
W_0^{1,p}(\mathcal{O})$. Therefore, taking $u^-$ as a test function,
we have
\[
\int_{\mathcal{O}}|\n u|^{p-2} \;\lan\n  u, \n u^-  \ran\,dx=
\lambda_1 \int_{\mathcal{O}}|u|^{p-2} \, u u^-\, dx.
\]
From this we obtain
\[
\int_{\mathcal{O}}|\n u^-|^{p} \,dx= \lambda_1
\int_{\mathcal{O}}|u^-|^{p}\, dx. 
\]
We cannot have $u^-\neq0$, for otherwise, from the variational
characterization of the first eigenvalue we can conclude that
$\lambda_1(\mathcal{O}) \le \lambda_1=\lambda_1(\Om)$. However, this
cannot happen,  $\mathcal{O}$ being a proper open subset of $\Om$ we
must have $\lambda_1(\Om) < \lambda_1(\mathcal{O})$.\\

\noindent \textbf{Existence.}\quad We denote by $f^\prime$ again a
$W^{1,p}(\mathcal{O})$ function whose trace on $\partial
\mathcal{O}$ is $f^\prime$. We can then obtain a weak solution of 
\eqref{nlpde} by minimizing the functional $J(w)=  
\int_{\mathcal{O}} |\n w|^p\, dx - \lambda_1(\Om)
\int_{\mathcal{O}}|w|^p\, dx$ on the affine space
$A:=W_0^{1,p}(\mathcal{O}) + f^\prime$. Indeed, if $w$  is  a
minimizer of $J$ then we shall have 
\be\label{varfrmltnnlpde}
\left. 0=\dfrac{d}{dt}\right\arrowvert_{t=0} J(w+ t \,
\phi)=\int_{\mathcal{O}}|\n w|^{p-2}\lan\n w, \n \vphi\ran\, dx -
\lambda_1\int_{\mathcal{O}}|w|^{p-2}w\, \vphi\, dx ~~\forall \; \vphi
\in \C_0^1(\mathcal{O}).
\ee
which is just the weak formulation of  \eqref{nlpde}. As $A$ is a
closed convex subset of the reflexive Banach space
$W^{1,p}(\mathcal{O})$,  for showing the existence of a minimizer  of
$J$ on $A$, it is enough to prove that $J$ is coercive and weakly
sequentially lower semi-continuous on $A$. \\ 

\noindent $J$ is weakly sequentially lower semi-continuous on
$A$: This is true since $\intom |\n w|^p dx$ is lower semicontinuous
for the weak topology on  $W^{1,p}(\mathcal{O})$ and $\intom |w|^p dx$
is continuous for the weak topology on  $W^{1,p}(\mathcal{O})$ due to
the compact inclusion of $W^{1,p}(\mathcal{O})$ in
$L^p(\mathcal{O})$. \\  
\noindent $J$ is coercive on $A$: \quad  Let $w_n := f^\prime + \vphi_n
\in A$ be a sequence such that $\|w_n\|_{W^{1,p}(\mathcal{O})}
\longrightarrow \infty$ as $n \to \infty$. If $\int_{\mathcal{O}}
|w_n|^p \,dx$ is  a bounded sequence, then the coercivity is immediate.

So, let us assume that $\int_{\mathcal{O}} |w_n|^p \,dx \to \infty$ as
$n \to \infty$.  We may write $w_n := f^\prime + \vphi_n$ with
$\vphi \in W^{1,p}_0(\mathcal{O})$. Let  $B_n :=
\dfrac{\int_{\mathcal{O}} |w_n|^p \,dx}{\int_{\mathcal{O}} |\phi_n|^p
  \,dx}$. It can be argued, using the 
triangle inequality, that  $\int_{\mathcal{O}} |\vphi_n|^p \,dx \to
\infty$ and $B_n \to 1$ as $n \to \infty$. From the Poincar\'e
inequality on $\mathcal{O}$, we conclude that 
$\int_{\mathcal{O}} |\n \vphi_n|^p \,dx \to \infty$ as $n \to
\infty$. Setting $A_n :=  \dfrac{\int_{\mathcal{O}} |\n
  w_n|^p\,dx}{\int_{\mathcal{O}} |\n \vphi_n|^p \,dx}$, we obtain 
using the triangle inequality, that $\int_{\mathcal{O}} |\n w_n|^p
\,dx \to \infty$ and $A_n \to 1$ as $n \to 
\infty$.
Now,  
\bea
\nonumber J(w_n) & = & \ds{A_n \, \lp \int_{\mathcal{O}}|\n \vphi_n|^p\,dx -
\lambda_1(\Omega) \dfrac{B_n}{A_n} \int_{\mathcal{O}}|\phi_n|^p\,dx \rp }\\
\label{coerJ1} &\ge & \ds{A_n \lp 1- \dfrac{B_n}{A_n}
  \dfrac{\lambda_1(\Om)}{\lambda_1(\mathcal{O})} \rp
  \int_{\mathcal{O}}|\n \vphi_n|^p\,dx  }  
\eea
where the last inequality has been obtained by applying Poincar\'e
inequality in the domain $\mathcal{O}$. Since we have $0 <
\lambda_1(\Omega)< \lambda_1(\mathcal{O})$, since $A_n$ and $B_n$
converge to $1$ as $n \to \infty$, it follows that $\ds{A_n \lp 
  1- \dfrac{B_n}{A_n} \dfrac{\lambda_1(\Om)}{\lambda_1(\mathcal{O})}
  \rp }$ is bounded below by a positive constant $C >0$. Once again,
we have the coercivity of $J$. \\ 

\noindent \textbf{Uniqueness.}\quad Suppose $u,v$ are two different
solutions of \eqref{varfrmltnnlpde} in $A$. Let $w_1:= \n \log u$ and
$w_2:=\n \log v$. As $f(x)=|x|^p$ is a strictly convex function we have  
\be\label{sbdffineq}
|w_1|^p \geq |w_2|^p+ p\, |w_2|^{p-2}\lan w_2,w_2-w_1\ran
\ee
and equality holds if and only if $w_1=w_2$. If we prove that
$w_1=w_2$ then we are done because in that case we will have $0=\n
\log u-\n \log v= \n \log\lp\frac{u}{v}\rp$. That is,
$\log\lp\frac{u}{v}\rp=k$ for some constant $k$. As a result we get $u=
e^k \, v$. But as $u \equiv v = f^\prime \not \equiv 0$ on $\partial
\mathcal{O}$ we get $u \equiv v$ in $\mathcal{O}$. Therefore, it suffices to prove that 
\be\label{sbdffeq}
|w_1|^p = |w_2|^p+ p\, |w_2|^{p-2}\lan w_2,w_2-w_1\ran.
\ee
The proof of (\ref{sbdffeq}) is the same as the proof of Lemma 3.1 in
Lindqvist \cite{Lindqvist}. We include the proof here for
completeness.  The function $u$ solves \eqref{nlpde}. We use $u-v^p \,
u^{1-p}$ as a test function in the equation for $u$. Similarly, we use
$v-u^p \, v^{1-p}$ as a test function in \eqref{nlpde} with $v$ as a
solution. Then we integrate by parts and sum the two identities. This
new identity can be reduced to  
\be\label{newidentity}
\begin{aligned}
0
& = & \int_{\mathcal{O}} \lc u^p \{|w_1|^p-|w_2|^p-p \, |w_2|^{p-2}\lan
w_2,w_2-w_1\ran\}\right.\\ 
&~   &+\left. v^p \{|w_2|^p-|w_1|^p-p \, |w_1|^{p-2}\lan
w_1,w_1-w_2\ran\}\rc\,dx.  
\end{aligned}
\ee
by using the following: 
\[
\n (u-v^p \, u^{1-p}) = \left\{1+(p-1) \lp\dfrac{v}{u}\rp^p\right\}\n
u-p\, \lp\dfrac{v}{u}\rp^{p-1}\n v,
\]
and, 
\[
\n (v-u^p \, v^{1-p}) = \left\{1+(p-1) \lp\dfrac{u}{v}\rp^p\right\}\n
v-p\, \lp\dfrac{u}{v}\rp^{p-1}\n u.
\] 
But by (\ref{sbdffineq}) the integrand in \eqref{newidentity} is
non-negative (being the sum of two non-negative terms) and so, it
follows from \eqref{newidentity} that this integrand is equal to zero
almost evrywhere in $\mathcal{O}$. Therefore, each of the terms in the
integrand must be zero. This proves (\ref{sbdffeq}). \hfill $\square$ 

\section{Shape derivatives of torsional rigidity and eigenvalue
  functionals} \label{shpdrvtvs}  
\setcounter{theorem}{0}
\setcounter{equation}{0}
Let $\Omega$ be a smooth domain in $\R^N$ and let $\mathcal{D}$ be a
domain such that $\Omega_t \subset \mathcal{D}$, for $t$ sufficiently
small, for the smooth perturbations $\Phi_t$ associated to a smooth
vector field $V$. Consider the Dirichlet boundary value
problem\label{DBVP4} on $\Om_t$ : 
\be\label{pdevardom}
\left.
\begin{aligned}
-\La_p u & =  1 ~& \mbox{ in } & ~\Om_t,  \\
u & =  0 ~& \mbox{ on } & ~\partial \Om_t. 
\end{aligned} 
~~~ \right\} 
\ee
Let $y_t \in \C^{1,\alpha}(\overline{\Om_t})$ be the unique solution
of problem (\ref{pdevardom}).  Throughout this section $y=y(\Om)$
denotes the unique solution of (\ref{pdevardom}) for $t=0$. Denote
$\left( y_t\circ \Phi_t \right) \arrowvert_\Om$ by $y^t$ ($t \in
\R$). We also denote the torsional rigidity $E(\Omega_t)$ by $E(t)$. 

\begin{propo} \label{exstshpdrvtvtrsnlrgdty} The shape derivative of
  the torsional rigidity functional $E(\Omega_t)$ exists at $t=0$ and  
\be\label{shpdrvtvtrsnrgdty}
\left. \dfrac{d}{dt}\right\arrowvert_{t=0} E(\Om_t)=\int_{\partial
  \Om}\left|\frac{\partial y}{\partial n} \right|^p\,\lan V, n\ran\; dS.
\ee
(Here, $n$ denotes unit outward normal on $\partial \Om$.)
\end{propo}
\noindent {\it Proof.}\quad Let $y$ be the unique solution of
\eqref{pdevardom} on $\Omega$ corresponding to $t=0$. 

\smallskip
\noi
{\sc \underline{Step 1}}: We first show that $y^t \longrightarrow y$
strongly in  $W_0^{1,p}(\Omega)$.  

This can be obtained using the $\Gamma$-convergence (cf. Attouch
\cite{At}, Braides \cite{Br}, Dal Maso \cite{DM}) of a suitable family
of functionals. Consider the following family of functionals
defined over $W_0^{1,p}(\Omega)$:   
\be\label{flyfnl}
\ds{F(t,y):= \frac{1}{p} \intom |B_t^*(\n y)|^p \gamma_t \, dx - \intom
   y \gamma_t \;  dx }
\ee
Since $B_t^*$ converges uniformly to $I$ and $\gamma_t$ converges
uniformly to the constant $1$, it is classical to show the
$\Gamma$-convergence of the family of convex integral functionals
$F(t, \cdot)$, as $t \to 0$, to the following functional  
\be\label{ltfnl}
\ds{F(y):= F(0,y)=\frac{1}{p} \intom |\n y|^p  \, dx - \intom \, y  \;
  dx } 
\ee
See Theorem 5.14 in Dal Maso \cite{DM} for instance. Furthermore, the 
family $F(t,\cdot)$ is equicoercive following the inequalities
\eqref{bdBt} and \eqref{asygammat}.  By standard results on
$\Gamma$-convergence (cf. Theorems 7.8 and 7.12 Dal Maso \cite{DM}),
the minimizer of $F(t,\cdot)$ converges weakly in $W^{1,p}_0(\Omega)$
to the minimizer of $F(\cdot)$ and the minima converge. Now, for each
$t \in \R$, $y_t$ satisfies the equation       
\be \label{vareqyt}
\intot |\n y_t|^{p-2}\lan\n y_t,\n \psi\ran\; dx = \intot 
\;\psi\; dx ~~~~~~ \forall \,\psi \in \C_0^\infty(\Om_t).
\ee
By the change of variable $\Phi_t : \Om \longrightarrow \Om_t$,
the equation (\ref{vareqyt}) can be re-written as  
\be \label{vareqyupt}
\intom |B_t^*(\n y^t)|^{p-2}\;\left<A_t(\n y^t)\,,\,\n \vphi
\right>\; dx=\intom \gamma_t\;\phi \; dx ~~\forall \; \vphi \in
\C_0^\infty(\Om).  
\ee
Therefore, $y^t$ satisfies :
\be\label{yupt}
\left. 
\begin{aligned}
-\mathrm{div}\lp|B_t^*(\n y^t)|^{p-2}A_t \lp\n y^t\rp\rp - \;\gamma_t
& =  0 ~& \mbox{ in } & ~\Om,  \\[1mm]  
y^t & =  0 ~& \mbox{ on } &  \partial \Om. 
\end{aligned} ~~~ \right\}
\ee
which is the Euler equation for the minimization of the convex
functional $F(t,\cdot)$ and therefore, $y^t$ is the minimizer of
$F(t,\cdot)$. Whereas,  $y$, being the solution of problem (\ref{pdevardom}) for
$t=0$, is the minimizer of $F$. So, by the $\Gamma$-convergence result,
we have the convergence of the minimum values
\be\label{cvmin}
\ds{\mathrm{lim}_{t \to 0} \lp \frac{1}{p} -1 \rp \intom |B_t^*(\n
  y^t)|^p \gamma_t \, dx = \lp \frac{1}{p} -1 \rp \intom | \n
  y|^p  \, dx \,.} 
\ee
and the weak convergence in $W^{1,p}_0(\Omega)$, as $t \to 0$ of $y_t$
to $y$. It remains to show the strong convergence. 

Since $B_t^*$ and $\gamma_t$ converge uniformly to $I$ and $1$
respectively, and $y_t$ remains bounded in $W^{1,p}_0(\Omega)$,
we can conclude from \eqref{cvmin} that
\be\label{strcvyt}
\ds{\mathrm{lim}_{t \to 0} \intom |\n y^t|^p \, dx =  \intom | \n 
  y|^p  \, dx \,.} 
\ee 
Therefore, since the $L^p$ norm is uniformly convex, we can conclude
from the weak convergence of $\n y^t$ to $\n y$ in $L^p(\Omega)$ and
the convergence of their norms (\ref{strcvyt}) that the convergence of
$\n y^t$ to $\n y$ is strong in $L^p(\Omega)$. By Poincar\'e
inequality, as the $L^p$ norm of the gradients is an equivalent norm
on $W^{1,p}_0(\Omega)$, we obtain the desired conclusion. 

\smallskip
\noi
{\sc \underline{Step 2}}: We observe that the torsional rigidity
$E(t)$ of the domain $\Omega_t$ is given by
\be\label{varchartrsnlrgdty}
E(t) =  \dfrac{p}{p-1} \sup_{\phi \in W^{1,p}_0(\Omega)} \left\{ \intom
\phi \gamma_t \;  dx - \frac{1}{p} \intom |B_t^*(\n \vphi)|^p \gamma_t
\, dx  \right\} =   \dfrac{p}{p-1} \sup_{\phi \in W^{1,p}_0(\Omega)}
(-F(t,\phi)) 
\ee
and the supremum is attained at $\phi = y^t$ for $y^t = y_t \circ
\Phi(t)$ and $y_t$ is the solution of \eqref{pdevardom} on $\Omega_t$.  

Indeed, the supremum in the above corresponds to the negative of the
infimum  in the following
\begin{eqnarray*}
&&\ds{\inf_{\phi \in W^{1,p}_0(\Omega)} \left\{ \frac{1}{p} \intom
    |B_t^*(\n \vphi)|^p \gamma_t \, dx  - \intom \vphi \gamma_t\;  dx 
\right\} }\\
&&=\ds{\inf_{\phi \in W^{1,p}_0(\Omega_t)} \left\{ \frac{1}{p}
    \int_{\Omega_t} |\n \vphi|^p  \, dx  - \int_{\Omega_t} \vphi \;  dx
    \right\} } 
\end{eqnarray*} 
and this is attained by $y_t$ which is the solution of the
Euler-Lagrange equation \eqref{pdevardom} on $\Omega_t$. We can
calculate this value which turns out to be
$$
\ds{F(t,y^t)= -\dfrac{p-1}{p} \int_{\Omega_t} y_t \ dx =
  -\dfrac{p-1}{p} E(t)}\,. 
$$
This proves our affirmation.

\smallskip
\noi
{\sc \underline{Step 3}}: We now show that the shape derivative
exists, that is the limit, $\ds{\mathrm{lim}_{t \to 0}
  \dfrac{E(t)-E(0)}{t}}$, exists and 
\be
\ds{\mathrm{lim}_{t \to 0} \dfrac{E(t)-E(0)}{t} = -\dfrac{p}{p-1}
  \dfrac{\partial F}{\partial t}(0,y)}
\ee 
We obtain from the variational characterization
\eqref{varchartrsnlrgdty} of $E(t)$ that  
\be \label{varcharimplctn1}
E(t)-E(0)\geq \dfrac{p}{p-1}\lp F(0,y) - F(t,y)\rp\,. 
\ee 
Thus, 
\be\label{part1}
\liminf_{t \downarrow 0} \dfrac{E(t)-E(0)}{t} \geq \dfrac{p}{p-1}
\lim_{t \downarrow 0} \dfrac{F(0,y)-F(t, y)}{t} = - \dfrac{p}{p-1}
\dfrac{\partial F}{\partial t}(0, y)\,. 
\ee
Once again by applying the variational characterization of
$E(t)$ we have
$$
E(t)-E(0) \leq \dfrac{p}{p-1} \lp F(0,y^t)-F(t, y^t) \rp \,. 
$$
Therefore, by applying the integral form of the mean value theorem in
the above in the first variable
\begin{eqnarray*}
\dfrac{E(t)-E(0)}{t} & \leq &  \dfrac{p}{p-1} \dfrac{ F(0,y^t) -
  F(t,y^t) }{t} \\  
& = & - \dfrac{p}{p-1} \int_0^1 \dfrac{\partial F(s \, t,
  y^t)}{\partial s} \ ds   
\end{eqnarray*}
In order to conclude the reverse inequality
\be \label{part2}
\limsup_{t \downarrow 0} \dfrac{E(t)-E(0)}{t} \leq
- \dfrac{p}{p-1} \dfrac{\partial F}{\partial t}(0, y) 
\ee
it is enough to show that
\be \label{claim-lt}    
\liminf_{t \downarrow 0} \dfrac{\partial F}{\partial t}(st, y^t)
= \dfrac{\partial F}{\partial t}(0, y) \mbox{ for every } s \in
     [0,1]\,.
\ee
By a straightforward computation it is seen that, for any $\psi \in
W^{1,p}_0(\Omega)$, we have 
\begin{eqnarray*}
&&\dfrac{\partial F}{\partial s}(s , \psi)  = 
\dfrac{\partial}{\partial s} \left( \frac{1}{p} \intom |B_s^*(\n
\psi)|^p \gamma_s \, dx - \intom \psi \gamma_s \;  dx \right)\\  
&=& \left(\int_\Om \dfrac{1}{p} \gamma'_s \, |B_s^* \n \psi|^p
  \ dx  + \int_\Om |\n \psi|^{p-2} \lan (B_s^*)'  \n \psi, B_s^* \n  
  \psi\ran\,dx -  \int_\Om \gamma'_s  \psi \ dx\right)\,. 
\end{eqnarray*}
So, in particular, by taking $\psi=y^t$ we get
\begin{eqnarray*}
\dfrac{\partial F}{\partial s}(st , y^t) =  \left(\int_\Om
\dfrac{1}{p} \gamma'_{st} \, |B_{st}^* \n y^t|^p \ dx  + \int_\Om |\n
y^t|^{p-2} \lan (B_{st}^*)'  \n y^t, B_{st}^* \n  y^t \ran\,dx -
\int_\Om \gamma'_{st}  y \ dx\right)\,.  
\end{eqnarray*}
Due to the strong convergence of $\n y^t$ to $\n y$ in $L^p(\Om)$
and the analyticity of $B_s$, $\gamma_s$ in $s$, it is now
straightforward to pass to the limit as $t \to 0$ and we obtain
easily, using \eqref{drvtvgammaat0} and \eqref{drvtvBt0},  
that 
\begin{eqnarray*}
\liminf_{t \downarrow 0} \dfrac{\partial F}{\partial s}(s t, y^t)
&=&  \left(\int_\Om \dfrac{1}{p} \mathrm{div} V  \, |\n y|^p \ dx  +
\int_\Om |\n y|^{p-2} \lan - (DV)^*  \n y\,,\, \n  y \ran\,dx -  \int_\Om  
\mathrm{div} V \,  y \ dx\right)\\
&= &   \dfrac{\partial F}{\partial t}(0, y)
\end{eqnarray*}
for every $s \in [0,1]$, proving the claim \eqref{claim-lt}.

\smallskip
\noi
\underline{{\sc Step 4}}: To obtain the expression for the shape
derivative \eqref{shpdrvtvtrsnrgdty}, it is enough to integrate by
parts in the term $ -  \int_\Om   \mathrm{div} V \,  y \ dx$ which
appears in the expression for $\dfrac{\partial F}{\partial t}(0,y)
= -\dfrac{p-1}{p}E'(0)$ and apply the Pucci-Serrin  identity
\eqref{PucciSerrinC1}.  \hfill $\square$

\bigskip
\noindent
We now recall the shape derivative for the eigenvalue
functional. Consider the eigenvalue problem: 
\be\label{evpvardom}
\left.
\begin{aligned}
-\La_p u & = \lambda |u|^{p-2}u ~& \mbox{ in } &~\Om_t , \\
u & = 0 ~& \mbox{ on } &~\partial \Om_t.
\end{aligned}  ~~~ \right\} 
\ee
The first eigenvalue $\lambda_1(t):= \lambda_1(\Om_t)$ is simple and
is characterized as the minimum of the problem 
\be\label{varchar1stevplapvardom}
\lambda_1(\Omega_t) :=
\inf \left\{ \left.\dfrac{\|\nabla \vphi\|^p_{L^p(\Om_t)}}{\|\phi
  \|^p_{L^p(\Om_t)}}\;\right|\; \vphi \in W^{1,p}_0(\Om_t)\right\}\,.
\ee
We fix $y_{1,t}:=y_1(\Om_t)$ to be a corresponding eigenfunction
which is positive (using the Krein-Rutman theorem) and normalize it to
satisfy
\be\label{nrmlztnt}
\int_{\Om_t} |y_{1,t}|^p \ dx = 1\,.
\ee
For $t=0$, we denote the corresponding eigenvalue and eigenfunction by
$\lambda_1$ and $y_1$ respectively.  

\begin{propo}\label{dffbltyev}
The map \,$t \longmapsto \lambda_1(t)$ is differentiable at
$\,0$ and 
\be\label{exprsndrvtvev}
\lambda_1^\prime(0)=-(p-1)\int_{\partial \Om}\left|\frac{\partial
  y_1}{\partial n} \right|^p\,\lan V, n\ran\; dS
\ee
\end{propo}
\noindent {\it Proof.} As we have mentioned before, this result has
been shown previously by de Lis and Garc\'ia-Meli\'an \cite{GMSlS} and
by Lamberti \cite{L} (also see Ly \cite{Ly}). This can also be proved
along the same lines as in Proposition
\ref{exstshpdrvtvtrsnlrgdty}. \hfill $\square$

\section{Main Results}\label{mainrslt}
\setcounter{theorem}{0}
\setcounter{equation}{0}
Let $0<r_0<r_1$, $B_1$ be the ball $B(0,r_1)$ and let $B_0$ be
any open ball of radius $r_0$ such that $\ds \overline{B_0}\subset
B_1$. Consider the family $\F = \left\{ B_1 \setminus
\overline{B_0}\right\}$ of domains in $\R^N$. We study the 
extrema of the functionals $E(\Om)$ and $\lambda_1(\Om)$ over $\F$,
associated to the problems (\ref{stationary}), (\ref{evp})
respectively. \\  
We state our main results :\\[3mm]
Put $\Om_0 = B(0,\,r_1) \setminus \overline{B(0, r_0)}$.
\begin{theorem}\label{trsnlrgdtyrslt}
The minimum value of the torsional rigidity functional $E(\Om)$ on
$\F$ is attained only when  $\Om = \Om_0$, i.e., when the balls are 
concentric. 
\end{theorem}
\begin{theorem}\label{evprslt}
The first Dirichlet eigenvalue $\lambda_1(\Om)$ is maximum on $\F$
only when  $\Om = \Om_0$, i.e., when the balls are concentric. 
\end{theorem}
Before proceeding to the proof we make the following observation and
reduction.  The functionals to be optimized are invariant under the
isometries of $\R^N$. Therefore, it is enough to study these
optimization problems for the class of domains
$\Om(s):=B_1\setminus \overline{B(s\,e_1,r_0)}$, $0 \leq s< r_1-r_0$
where $e_1$ is the unit vector in the direction of the first
coordinate axis. In order to study the optimality of the domain
$\Om(s)$ in the class $\F$ we need to study perturbations of the
domain which correspond to translations of the inner ball along the
direction of the first coordinate axis. For this purpose we consider a
smooth vector field $V(x)= \rho(x)\, e_1 \; \forall \; x 
\in B_1$ where $\rho: \R^N \rightarrow [0,1]$ is a smooth function with compact support in
$B_1$ such that $\rho \equiv 1$ on a neighborhood of $B(s\, e_1,
r_0)$. Let $\{\Phi_t\}_{t \in \R}$ be the one-parameter family of 
diffeomorphisms of $B_1$ associated with $V$. We see that, for $t$
sufficiently close to $0$, $\Phi_t(\Om(s)))= \Om(s+t)$. So, if we
define $j,j_1: (r_0-r_1,r_1-r_0) \rightarrow \R$ as follows: 
\[
j(s)=E(\Om(s)) ~~~~\mbox{ and }~~~~j_1(s)=\lambda_1(\Om(s))\,
\] 
we see that the minimization of $E$ in the class $\F$ corresponds to
studying the minimum of $j$ on the interval $(r_0-r_1,r_1-r_0)$
and that the problem of maximization of $\lambda_1$ in the class $\F$
corresponds to studying the maximum of $j_1$ on the interval
$(r_0-r_1,r_1-r_0)$. Also, the shape derivative of $E$ and $\lambda_1$
at $\Om(s)$ for the vector field $V$ are the ordinary derivatives at
$s$ of $j$ and $j_1$ respectively. We have seen in Proposition
\ref{exstshpdrvtvtrsnlrgdty} and \ref{dffbltyev} that these shape
derivatives exist and so the derivative of both $j$ and $j_1$
exist. The optimization problems can be studied by analyzing the sign
of the derivatives of $j$ and $j_1$. 

First, we note that both $j$ and $j_1$ are even functions and since
they are differentiable, we have $j^\prime(0)=0=j_1^\prime(0) $.   

We shall adopt the following notations. Given $s$ in $(0,r_1-r_0)$ we
simply denote $\Om(s)$ as $\Om$ and $B(s\,e_1, r_0)$ as $B_0$ and $n$
shall denote the unit outward normal to $\Om$ on $\partial \Om$. Let $H$ denote
the hyperplane $H:=\{x=(x_1,x_2, \ldots, x_N)\in \R^N \,|\,
x_1=s\}$. Let $r_H$ be reflection function about $H$. We define
$\mathcal{O}$ to be the subdomain $\mathcal{O}:=\{x \in \Om \,|\, x_1
>s\}$ in $\Om$. Then we see that the reflection of $\mathcal{O}$ about
$H$, namely $\mathcal{O'}=r_H(\mathcal{O})$ is contained in $B_1$, 
whereas $B_0$ is symmetric with respect to $H$.  Thus $\mathcal{O'} \subset \Om$. For $x \in \mathcal
{O}$, let $x^\prime$ denote the reflection of $x$ about $H$, namely,
the point $r_H(x)$. With these notations, if $y$ be the solution of the
equation \eqref{stationary} in $\Omega$ then, from the expression of 
the shape derivative \eqref{shpdrvtvtrsnrgdty} for $E$  we obtain that  
\be\label{drvtvj} 
j^\prime (s) = \int_{x \in \partial B_0}
\left|\frac{\partial y}{\partial n}(x) \right|^p \; n_1(x) \; dS 
\ee
since $V$ is zero on $\partial B_1$ and since, for all $x \in
\partial B_0 $, $<V,n>(x)=\rho e_1 \cdot n(x) 
=n_1(x)$, the first component of the normal vector. Similarly, if
$y_1$ be the solution of \eqref{evp} in $\Omega$ then, from the
expression of the shape derivative of
$\lambda_1$, viz. \eqref{exprsndrvtvev}, we obtain that 
\be\label{drvtvj1} 
j_1^\prime (s)  =  -(p-1)\int_{x \in \partial B_0}
\left|\frac{\partial y_1}{\partial n}(x)\right|^p \; n_1(x) \; dS
\ee 

\noindent {\bf Proof of Theorem \ref{trsnlrgdtyrslt}} Let $y$ be the
solution of the boundary value problem \eqref{stationary} in
$\Omega$. We recall that $y \in \C^{1, \alpha} (\ov{\Omega})$ by
regularity results in Tolksdorff \cite{Tolk2}, and by the strong
maximum principle (cf. Theorem 5, Vazquez \cite{Vaz}) we have $y>0$ in
$\Omega$. We now consider the subdomain $\mathcal{O}$ and let us
define $\tilde{y}$ on $\mathcal{O}$ by $\tilde{y}(x):= y(x^\prime)$
the value of $y$ at the reflection $x'$ of $x$ about $H$. Let us note
that $\frac{\partial \tilde{y}}{\partial n}(x) = \frac{\partial
  y}{\partial n}(x^\prime)$ and $n_1(x^\prime)=-n_1(x)$ for all $x \in
\partial B_0$. Now, we may rewrite the expression
\eqref{drvtvj} as follows: 
\bea
\nonumber j^\prime (s) & = & \int_{x \in \partial B_0\cap \partial \mathcal {O}}
\left\{\left|\frac{\partial y}{\partial n}(x) 
\right|^p n_1(x) \; dS + \int_{x' \in \partial B_0\cap \partial \mathcal {O'}}
\left|\frac{\partial y}{\partial n}(x^\prime) \right|^p
\right\}\; n_1(x^\prime) \; dS  \\
\label{exprsndrvtvj}& = & \int_{x \in \partial B_0\cap \partial \mathcal {O}}
\left\{\left|\frac{\partial y}{\partial n}(x) 
\right|^p-\left|\frac{\partial \tilde{y}}{\partial n}(x) \right|^p
\right\}\; n_1(x) \; dS . 
\eea
We shall show that $j^\prime(s) \geq 0 \; \forall \; s \in [0,
  r_1-r_0)$ and is zero only if $s=0$. We have already observed that
  $j^\prime(0)=0$ by symmetry considerations. It is clear that $n_1(x)<0\; \forall \; x \in \partial \mathcal{O} \cap \partial B_0 \cap H^c$. So when $s \neq 0$, we
  shall prove that $j^\prime(s)>0$ by showing that   
\be\label{inequality}
\frac{\partial \tilde{y}}{\partial n}(x) < \frac{\partial y}{\partial
  n}(x)<0 ~~~\, \forall\; x \in \partial \mathcal{O} \cap \partial B_0 \cap H^c.
\ee 
We shall prove inequality (\ref{inequality}) in a few steps. \\ 
\noindent {\underline{\sc Step 1}: } First we prove that $\dfrac{\partial
  y}{\partial n}<0$ on $\partial B_0$. \\ 
We begin by noticing that at every point $x_0$ on $\partial B_0$, the
interior sphere property holds, that is, there exists an open ball 
$B= B_R(z_0) \subset \Omega$ such that $\partial B \cap \partial B_0 = 
\{x_0\}$ and the unit outward normal $n$ to $\Omega$ and to $B$ coincide
at $x_0$. For $K >0$ and $\alpha >0$ we define a function $b : B
\rightarrow \R$ as $b(x) = K \lp e^{-\alpha |x-z_0|^2} - e^{-\alpha
  R^2} \rp$. We have $b >0$ in $B_R(z_0) \setminus
B_{\frac{R}{2}}(z_0)$, $b(x_0) =0$, in fact, $b \equiv 0$ on $\partial
B$ and that $\dfrac{\partial b}{\partial n}(x_0)<0$. Moreover, it can
be shown that $$-\La_p b(x) =- 2K e^{-\alpha |x-z_0|^2} |\n b|^{p-2}
\lp 2 \alpha^2 |x-z_0|^2 -N \, \alpha \rp. $$  
We may therefore, choose $\alpha$ large enough (independent of $K$) so that 
\be\label{Lapb} 
- \La_p b \leq 0 \mbox{ in } B_R(z_0) \setminus B_{\frac{R}{2}}(z_0).
\ee 
We know that $y$ satisfies (\ref{stationary}).  Since $y$ is bounded
below by a positive constant on $\partial B_{\frac{R}{2}}(z_0)$, we
may choose $K$ small enough so that $b \leq y$ on $\partial
B_{\frac{R}{2}}(z_0)$. Thus we have  
\be\label{bleqy}
b \leq y \mbox{ on } \partial \lp B_R(z_0) \setminus
B_{\frac{R}{2}}(z_0) \rp,~ -\La_p b \leq 0 ~\mbox{ and } -\La_p y >0
\mbox{ in } B_R(z_0) \setminus B_{\frac{R}{2}}(z_0).  
\ee
Then by the WCP of Tolksdorff \cite{Tolksdorff} we have $b \leq y$ in
$B_R(z_0) \setminus B_{\frac{R}{2}}(z_0)$.  This, along with $b(x_0)=0
=y(x_0)$, implies that $\dfrac{\partial b}{\partial n}(x_0) \geq
\dfrac{\partial y}{\partial n}(x_0)$. Since $\dfrac{\partial
  b}{\partial n}(x_0) <0$ we get $\dfrac{\partial y}{\partial
  n}(x_0)<0$. We have thus proved that $\dfrac{\partial y}{\partial
  n}<0$ on $\partial B_0 \cap \partial \mathcal{O} \cap H^c$. 

\noindent {\underline{\sc Step 2}: } Now we prove the first inequality in
(\ref{inequality}). \\ 
On $\mathcal{O}$, the function $y$ satisfies 
\begin{equation*}
\begin{aligned}
-\La_p y& =1 & \mbox{ in } & \mathcal{O},\\
y & =0 & \mbox{ on } & \partial \mathcal{O} \cap \partial B_0
,\\
y & =y^* & \mbox{ on } & \partial \mathcal{O} \cap H,\\
y & =0 & \mbox{ on } & \partial \mathcal{O} \cap \partial B_1;
\end{aligned}
\end{equation*}
while $\tilde{y}$ satisfies   
\begin{equation*}
\begin{aligned}
-\La_p \tilde{y}& =1 & \mbox{ in } & \mathcal{O},\\
\tilde{y} & =0 & \mbox{ on } & \partial \mathcal{O} \cap \partial B_0
,\\
\tilde{y} & =y^* & \mbox{ on } & \partial \mathcal{O} \cap H,\\
\tilde{y} & >0 & \mbox{ on } & \partial \mathcal{O} \cap \partial B_1;
\end{aligned}
\end{equation*}
where $y^*$ denote the common value of $y$ and $\tilde{y}$ on
$H$.
Therefore,  we have
\begin{equation*}
\begin{aligned}
-\La_p \tilde{y}& =-\La_p y & \mbox{ in } & \mathcal{O},\\
\tilde{y} & \geq y & \mbox{ on } & \partial \mathcal{O}.
\end{aligned}
\end{equation*}
So, by the WCP  of Tolksdorff \cite{Tolksdorff} we get, 
$\tilde{y} \geq y$ in $\mathcal{O}$. Since, $\tilde{y}\equiv 0 \equiv
y$ on $\partial B_0 \cap \partial \mathcal{O}$ we conclude that 
\be\label{wkineqy-1}
\frac{\partial \tilde{y}}{\partial n}(x) \le \frac{\partial
  y}{\partial n}(x)  \mbox{ for all } x \in \partial B_0 \cap \partial
\mathcal{O}\,. 
\ee
By the result of Step 1 and (\ref{wkineqy-1}), we can obtain a
neighborhood $\mathcal{N}$ of $\partial B_0 \cap \partial 
\mathcal{O}$ and positive numbers $\eta, \epsilon_0$ small enough so that 
\be\label{lwbdcvxcmbntngr}
|t \n \tilde{y}(x)+(1- t) \n y(x)| \geq \eta, ~ t \tilde{y}(x)+
(1-t)  y(x) \leq \epsilon_0 ~ \forall \;t \in [0,1], \; \forall \; x
\in \overline{\mathcal{N}}.
\ee
Let $w :=\tilde{y}-y $, then $w \geq 0$ on $\partial \mathcal{O}$ with
$w =0$ on $\partial \mathcal{O} \cap \lp H \cup \partial B_0 \rp$ and
$w >0$ on $\partial \mathcal{O} \cap \partial B_1$.  
Let $A : \R^N \rightarrow \R^N$ be the map as defined in
\eqref{2.10}. Then, we have   
$$div(A(\n \tilde{y})-A(\n y))=\La_p \tilde{y}- \La_p y= 0 \mbox{ in }
\mathcal{N}.$$ 
By Mean Value Theorem we get,
$$\sum_{i=1}^N \dfrac{\partial}{\partial x_i} \left< \int_0^1 (\n
A_i)(t \n \tilde{y} + (1-t) \n y)\, dt ,\n w \right> =0.$$ 
Let $a_{ij}(x) = \int_0^1 \frac{\partial A_i}{\partial x_j}(t \n
\tilde{y}(x) + (1-t) \n y(x))\, dt$. Then, $w$ satisfies  
\be 
\begin{aligned}
-\sum_{i,j =1}^N \dfrac{\partial}{\partial x_i}\lp a_{ij}
\dfrac{\partial w}{\partial x_j}\rp & =0 & \mbox{ in } &
\mathcal{N},\\ 
w & \geq 0 & \mbox{ in } & \mathcal{N},\\
w & = 0 & \mbox{ on } & \partial \mathcal{N} \cap \partial B_0.
\end{aligned} 
\ee 
By \eqref{2.11},  $\left[ a_{ij}(x)\right]_{i,j =1}^N \ge \min\{1,p-1\}
\int_0^1 |t \n \tilde{y}(x) + (1-t) \n y(x))|^{p-2}\, dt$  and so, by
\eqref{lwbdcvxcmbntngr} we conclude that $\left[ a_{ij}(x)\right]_{i,j
  =1}^N$ is a uniformly positive definite matrix when $x \in
\overline{\mathcal{N}}$. Therefore, by the maximum principle for
uniformly elliptic operators (cf. Theorem 5, Ch. 2, Protter and
Weinberger \cite{PW}), since $w$ is a non-constant function, it
follows that the minimum of $w$ will be attained on $\partial
\mathcal{N}$. Since $\inf_{\partial \mathcal{N}} w =0$ it follows that
$w > 0$ in $\mathcal{N}$. Further, by the same argument as in the
Hopf's Lemma for uniformly elliptic operators (cf. Theorem 7, Ch. 2,
Protter and Weinberger \cite{PW}) we have $\dfrac{\partial w}{\partial
  n} <0$ on $\partial \mathcal{N} \cap \partial B_0$. That is,
$\forall\; x \in \partial \mathcal{O} \cap \partial B_0$ 
we have the following:  
$$
\frac{\partial \tilde{y}}{\partial n}(x)<\frac{\partial y}{\partial
  n}(x)<0. 
$$
\hfill $\square$\\

\noindent {\bf Proof of Theorem \ref{evprslt}}  
Recall from Section
\ref{shpdrvtvs} that $y_1$ is the principal eigenfunction of
\eqref{evp} (that is, the unique solution of \eqref{evp} for $\lambda=\lambda_1(\Omega)$)
 characterized by  $y_1   > 0$ in $\Omega$ and $\int_\Omega y_1^p \, dx =1$. We now 
consider the subdomain $\mathcal{O} $ and let us define $\tilde{y}_1$
on $\mathcal{O}$ by $\tilde{y}_1(x):= 
y_1(x^\prime)$ the value of $y_1$ at the reflection $x'$ of $x$ about
$H$. Let us note that $\frac{\partial \tilde{y}_1}{\partial n}(x) =
\frac{\partial y_1}{\partial n}(x^\prime)$ and $n_1(x^\prime)=-n_1(x)$
for all $x \in \partial B_0$. Now, we may rewrite the
expression \eqref{drvtvj1} as follows: 
\bea
\nonumber j_1^\prime (s) & = & -(p-1) \int_{x \in \partial B_0\cap \partial \mathcal {O}}
\left\{\left|\frac{\partial y_1}{\partial n}(x) 
\right|^p n_1(x) \; dS -(p-1) \int_{x' \in \partial B_0\cap \partial \mathcal {O'}}
\left|\frac{\partial y_1}{\partial n}(x^\prime) \right|^p
\right\}\; n_1(x^\prime) \; dS  \\
\label{exprsndrvtvj1}& = & -(p-1) \int_{x \in \partial B_0\cap \partial \mathcal {O}}
\left\{\left|\frac{\partial y_1}{\partial n}(x) 
\right|^p-\left|\frac{\partial \tilde{y}_1}{\partial n}(x) \right|^p
\right\}\; n_1(x) \; dS . 
\eea
We shall show that $j_1^\prime(s) \leq 0 \; \forall \; s \in [0,
  r_1-r_0)$ and is zero only if $s=0$. We have already observed that
  $j_1^\prime(0)=0$ by symmetry considerations. It is clear that $n_1(x)<0\; \forall \; x \in \partial \mathcal{O} \cap \partial B_0 \cap H^c$. So when $s \neq 0$, we
  shall prove that $j_1^\prime(s)<0$ by showing that   
\be\label{inequality1}
\frac{\partial \tilde{y}_1}{\partial n}(x) < \frac{\partial
  y_1}{\partial n}(x)<0 ~~~\, \forall\; x \in \partial \mathcal{O}
\cap \partial B_0 \cap H^c. 
\ee 
We shall prove inequality (\ref{inequality1}) in a few steps. \\
\noindent {\underline{\sc Step 1}: } First we prove that $\dfrac{\partial
  y_1}{\partial n}<0$ on $\partial B_0$. \\ 
As in the proof of Theorem \ref{trsnlrgdtyrslt} we begin by noticing
that at every point $x_0$ on $\partial B_0$, the interior sphere
property holds, that is, 
there exists an open ball 
$B= B_R(z_0) \subset \Omega$ such that $\partial B \cap \partial B_0 = 
\{x_0\}$ and the unit outward normal $n$ to $\Omega$ and to $B$ coincide at $x_0$. 
We recall that $y_1$ satisfies (\ref{evp}) and $y_1 >0$ in
$\Omega$. We construct an auxiliary function $b$, as in the proof of
Theorem \ref{trsnlrgdtyrslt}, with $\alpha$ sufficiently large so that
$- \La_p b \leq 0$ on $B_R(z_0) \setminus B_{\frac{R}{2}}(z_0)$ and
$K$ small so that $b \leq y_1$ in $\partial \lp B_R(z_0) \setminus
B_{\frac{R}{2}}(z_0) \rp$, to obtain $b \leq y_1$ in $B_R(z_0)
\setminus B_{\frac{R}{2}}(z_0)$ and consequently to obtain
$\dfrac{\partial y_1}{\partial n}<0$ on $\partial B_0$.  

\noindent {\underline{\sc Step 2}: } Now we prove the first inequality
in (\ref{inequality1}).\\ 
On $\mathcal{O}$, the function $y_1$ satisfies 
\begin{equation*}
\begin{aligned}
-\La_p y_1& =\lambda_1\, y_1^{p-1} & \mbox{ in } &
\mathcal{O},\\ 
y_1 & =0 & \mbox{ on } & \partial \mathcal{O} \cap \partial
B_0,\\ 
y_1 & =y_1^* & \mbox{ on } & \partial \mathcal{O} \cap H,\\
y_1 & =0 & \mbox{ on } & \partial \mathcal{O} \cap \partial B_1;
\end{aligned}
\end{equation*}
whereas $\tilde{y}_1$ satisfies 
\begin{equation*}
\begin{aligned}
-\La_p \tilde{y}_1& =\lambda_1\, \tilde{y}_1^{p-1} & \mbox{ in } &
\mathcal{O},\\ 
\tilde{y}_1 & =0 & \mbox{ on } & \partial \mathcal{O} \cap \partial
B_0,\\ 
\tilde{y}_1 & =y_1^* & \mbox{ on } & \partial \mathcal{O} \cap H,\\
\tilde{y}_1 & >0 & \mbox{ on } & \partial \mathcal{O} \cap \partial B_1;
\end{aligned}
\end{equation*}
where $y_1^*$ denotes the common value of $y_1$ and $\tilde{y}_1$ on
$H$. Therefore, we have  
\begin{equation*}
\begin{aligned}
 -\La_p \tilde{y}_1& =\lambda_1\, \tilde{y}_1^{p-1} & \mbox{ in } &
 \mathcal{O},\\ 
\tilde{y}_1 & \geq y_1 & \mbox{ on } & \partial \mathcal{O}.
\end{aligned}
\end{equation*} 
Therefore, by the WCP (cf. Theorem \ref{WCPbetainc} proved in 
Section \ref{comparison-plap}) for $\frac{\partial \beta}{\partial u}
\geq 0$  we get, $\tilde{y}_1 \geq y_1$ in $\mathcal{O}$. 
 Since, $\tilde{y}_1 \equiv 0 \equiv y_1$ on $\partial B_0 \cap
 \partial \mathcal{O}$ we conclude that  
\be\label{wkineqy1-1}
\frac{\partial \tilde{y}_1}{\partial n}(x) \le \frac{\partial
  y_1}{\partial n}(x)  \mbox{ for all } x \in \partial B_0 \cap \partial
\mathcal{O}\,. 
\ee
By the result of Step 1 and (\ref{wkineqy1-1}), we can obtain a
neighborhood $\mathcal{N}$ of $\partial B_0 \cap \partial 
\mathcal{O}$ and positive numbers $\eta, \epsilon_0$ small enough so that 
$$|t \n \tilde{y}_1(x)+(1- t) \n y_1(x)| \geq \eta, ~ t
\tilde{y}_1(x)+ (1-t)  y_1(x) \leq \epsilon_0 ~ \forall \;t \in [0,1],
\; \forall \; x \in \overline{\mathcal{N}}.$$ 
Let $w :=\tilde{y}_1-y_ 1$, then $w \geq 0$ on $\partial \mathcal{O}$
with $w =0$ on $\partial \mathcal{O} \cap \lp H \cup \partial B_0 \rp$
and $w >0$ on $\partial \mathcal{O} \cap \partial B_1$. 
We have  
$$-div(A(\n \tilde{y}_1)-A(\n y_1))=-\La_p \tilde{y}_1+ \La_p y_1=
\lambda_1 \lp \tilde{y}_1^{p-1} -y_1^{p-1}\rp \geq 0 \mbox{ in }
\mathcal{N},$$ where the map $A$ is as defined in \eqref{2.10}. By Mean Value Theorem we get,
$$-\sum_{i=1}^N \dfrac{\partial}{\partial x_i} \left< \int_0^1 (\n
A_i)(t \n \tilde{y}_1 + (1-t) \n y_1)\, dt ,\n w \right> \geq 0.$$ 
Let $a_{ij}(x) = \int_0^1 \frac{\partial A_i}{\partial x_j}(t \n
\tilde{y}_1(x) + (1-t) \n y_1(x))\, dt$, then $w$ satisfies  
\be 
\begin{aligned}
-\sum_{i,j =1}^N \dfrac{\partial}{\partial x_i}\lp a_{ij}
\dfrac{\partial w}{\partial x_j}\rp & \geq 0 & \mbox{ in } &
\mathcal{N},\\ 
w & \geq 0 & \mbox{ in } & \mathcal{N},\\
w & = 0 & \mbox{ on } & \partial \mathcal{N} \cap \partial B_0.
\end{aligned} 
\ee 
As in the proof of Proposition \ref{trsnlrgdtyrslt}, we observe that
the matrix $\left[ a_{ij}(x)\right]_{i,j =1}^N$ is a uniformly
positive definite matrix when $x \in \overline{\mathcal{N}}$. Then by
the maximum principle for uniformly elliptic operators (cf. Theorem 5,
Ch. 2, Protter and Weinberger \cite{PW}), since $w$ is a non-constant
function, it follows that the minimum of $w$ will be attained on
$\partial \mathcal{N}$. Since $\inf_{\partial \mathcal{N}} w =0$ it
follows that $w > 0$ in $\mathcal{N}$. Further, by the same argument
as in the Hopf's Lemma for uniformly elliptic operators (cf. Theorem
7, Ch. 2, Protter and Weinberger \cite{PW}), we have $\dfrac{\partial
  w}{\partial n} <0$ on $\partial \mathcal{N} \cap \partial B_0$. That
is, $\forall\; x \in \partial \mathcal{O} \cap \partial B_0$ 
we have the following:  
$$
\frac{\partial \tilde{y}_1}{\partial n}(x)<\frac{\partial y_1}{\partial
  n}(x)<0. 
$$
\hfill $\square$
\begin{acknowledgements} The authors thank S. Kesavan (Institute of
  Mathematical Sciences, Chennai, India) for fruitful discussions. The
  authors also thank S. Prashanth (TIFR-CAM, Bangalore, India) for
  sharing his  notes on comparison principles.  
\end{acknowledgements}


\begin{thebibliography}{99} 
\bibitem{Anane} Anane A. , Simplicit\'{e} et isolation de la
  premi\`{e}re valeur propre du p-laplacien avec poids, 
    C. R. Acad. Sci. Paris S\'{e}r. I Math. 305 (16), 725--728 (1987). 

\bibitem{Anisa-Aithal} Anisa M.H.C. and Aithal A.R. , On two
  functionals connected to the Laplacian in a class of doubly
  connected domains in space-forms, Proc. Indian
    Acad. Sci. (Math. Sci.), 115(1), 93--102 (2005). 


\bibitem{At} Attouch H. , Variational Convergence for Functions and
Operators, Applicable Math. Series, Pitman, Boston, 1984

\bibitem{Berestycki-Nirenberg} Berestycki H. and Nirenberg L. , On the
  moving plane method and the sliding method, Boll. Soc. Brasiliera
  Mat. Nova Ser. , 22, 1--37 (1991).

\bibitem{Br} Braides A. , Gamma-Convergence for Beginners
Oxford Lecture Series in Mathematics and Its Applications, 22,
Clarendon Press, 2002.  

\bibitem{Cuesta} Cuesta M. and Tak\'ac P. , A strong comparison principle
  for positive solutions of degenerate elliptic equations, 
    Differential and Integral Equations, 13 (4-6), 721--746 (2000).  

\bibitem{DM} Dal Maso G. , An Introduction to $\Gamma$-convergence,
PNLDE 8, Birkh\"auser, 1993.

\bibitem{DMS} Degiovanni M. , Musesti A. and Squassina M. , On the
regularity of solutions in the Pucci-Serrin identity, 
  Calc. Var. , 18, 317--334 (2003).


\bibitem{dPFB} Del Pezzo L.M. and Fern\'andez-Bonder J. , Some
  optimization problems for the $p$-Laplacian type equations,
  Appl. Math. Optim. , 59, 365--381 (2009). 
  

\bibitem{GMSlS} Garc\'ia-Meli\'an J. and Sabina de Lis J. , On the
  perturbation of eigenvalues for the $p$-Laplacian,
  C.R. Acad. Sci. Paris, t. 332, 893--898 (2001), .  

\bibitem{Gidas-Ni-Nirenberg} Gidas B. , Ni W.M. and Nirenberg L. ,
  Symmetry and related properties via the maximum principle,
  Comm. Math. Phys., 68, 209-243 (1979). 


\bibitem{HKK} Harrell E.M. , Kr\"oger P. and Kurata K. , On the
  placement of an obstacle or a well as to optimize the fundamental
  eigenvalue, SIAM J. Math. Anal. , 33(1), 240-259 (2001). 

\bibitem{Kesavan} Kesavan S. , On two functionals connected to the
  Laplacian in a class of doubly connected domains, Proc.
    Roy. Soc. Edinburgh Sec. A., 133, 617--624 (2003). 

\bibitem{L} Lamberti P.D. , A differentiability result for the first
  eigenvalue of the $p$-Laplacian upon domain perturbation, Nonlinear
  analysis and applications: to V. Lakshmikantham on his  
80th birthday vol. 1, 2, Kluwer Acad. Publ., Dordrecht, 741--754
(2003). 
  

\bibitem{Lindqvist} Lindqvist P. , On the equation
  $\operatorname{\mathrm{div}}(|\nabla u|^{p-2}\nabla u)+ \lambda|u|^{p-2}u=0$,
  Proc. Amer. Math. Soc. , 109(1), 157--164 (1990). 

\bibitem{Ly} Ly I. , The first eigenvalue of the $p$-Laplacian
  operator, Journal Ineq. Pure Appl. Math. , 6(3), (2005).


\bibitem{PW} Protter M. and Weinberger H. , Maximum Principles in
  Differential Equations, Springer-Verlag New York, 1999. 

\bibitem{Ramm-Shivakumar} Ramm A.G. and Shivakumar P.N. ,
  Inequalities for the minimal eigenvalue of the Laplacian in an
  annulus, Math. Inequalities and Appl., 1(4), 559--563  (1998). 

\bibitem{S} Simon J. , Differentiation with respect to the domain in
  boundary value problems, Numer. Funct. Anal. and Optimiz. ,
  2(7--8), 649--687 (1980).

\bibitem{Sokolowski-Zolesio} Sokolowski J. and Zolesio J.P. ,
  Introduction to shape optimization: shape sensitivity analysis, 
  Springer series in computational mathematics, 10,
  Springer-Verlag, Berlin, New York, 1992. 

\bibitem{Tolksdorff} Tolksdorff P. , On the Dirichlet Problem for
  quasilinear equations, Comm. in PDEs, 8(7), 773-817 (1983). 

\bibitem{Tolk2} Tolksdorff P. , Regularity for a more general class
  of quasilinear elliptic equations, J. Differential
    Equations, 51, 126--150 (1984). 

\bibitem{Vaz} Vazquez J.L. , A strong maximum principle for some
  quasilinear elliptic equations, Appl. Math. Optim. , 12, 191--202
  (1984).     
\end{thebibliography}


\end{document}